\documentclass[leqno,12pt]{article}
\usepackage{latexsym}
\usepackage{hyperref}
\usepackage[pdftex]{graphicx}
\newcommand{\nc}{\newcommand}
\nc{\nt}{\newtheorem} \nt{thm}{Theorem}[section]
\nt{cor}[thm]{Corollary} \nt{prop}[thm]{Proposition}
\nt{lem}[thm]{Lemma} \nt{defn}[thm]{Definition}
\nt{exa}[thm]{Example} \nt{ass}[thm]{Assumption}
\nt{alg}[thm]{Algorithm} \nt{com}[thm]{Comment for Dennis}
\nt{coml}[thm]{Comment for Adrian} \nc{\ip}[2]{\mbox{$\langle #1,#2
\rangle$}} \nc{\pf}{\noindent{\bf Proof\ \ }}
\nc{\finpf}{\hfill{$\Box$}\linespace}
\nc{\linespace}{\vspace{\baselineskip} \noindent} \nc{\bx}{\bar x}
\nc{\R}{{\bf R}} \nc{\Rn}{{\bf R}^n} \nc{\Sn}{{\bf S}^n}
\nc{\Rm}{{\bf R}^m} \nc{\E}{{\bf E}} \nc{\e}{\epsilon}
\nc{\rt}{\rightarrow} \nc{\Diag}{{\rm Diag}\,} \nc{\tra}{\mbox{\rm
tr}\,}
\def\tto{\;{\lower 1pt \hbox{$\rightarrow$}}\kern -12pt
           \hbox{\raise 2.8pt \hbox{$\rightarrow$}}\;}
\newenvironment{myequation}{\setcounter{equation}{\value{thm}}
   \begin{equation}}{\addtocounter{thm}{1}\end{equation}}

\nc{\bmye}{\begin{myequation}} \nc{\emye}{\end{myequation}}
\newcounter{algorithm}
\renewcommand{\thealgorithm}{\arabic{section}.\arabic{algorithm}}
\nc{\newalgorithm}{\stepcounter{algorithm}\thealgorithm}

\begin{document}
\title{Randomized Methods for Linear Constraints:  Convergence Rates and Conditioning}
\author{
D. Leventhal\thanks{School of Operations Research and Information
Engineering, Cornell University, Ithaca, NY 14853, U.S.A.
\texttt{leventhal\char64 orie.cornell.edu }} \and A.S.
Lewis\thanks{School of Operations Research and Information
Engineering, Cornell University, Ithaca, NY 14853, U.S.A.
\texttt{aslewis\char64 orie.cornell.edu}}}

\maketitle

\noindent {\bf Key words:} coordinate descent, linear constraint, condition number, randomization, error bound, iterated projections, averaged projections, distance to ill-posedness, metric regularity
 \\
{\bf AMS 2000 Subject Classification: 15A12, 15A39, 65F10, 90C25}

\begin{abstract}
We study randomized variants of two classical algorithms:  coordinate descent for systems of linear equations and iterated projections for systems of linear inequalities.  Expanding on a recent randomized iterated projection algorithm of Strohmer and Vershynin for systems of linear equations,   we show that, under appropriate probability distributions, the linear rates of convergence (in expectation) can be bounded in terms of natural linear-algebraic condition numbers for the problems.  We relate these condition measures to distances to ill-posedness, and discuss generalizations to convex systems under metric regularity assumptions.
\end{abstract}

\section{Introduction}
The condition number of a problem measures the sensitivity of a
solution to small perturbations in its input data. For many problems
that arise in numerical analysis, there is often a simple
relationship between the condition number of a problem instance and
the distance to the set of ill-posed problems---those problem
instances whose condition numbers are infinite \cite{Demmel2}. For example, with
respect to the problem of inverting a matrix $A$, it is known (see
\cite{Horn}, for example) that if $A$ is perturbed to $A+E$ for
sufficiently small $E$, then
\[
\frac{\|(A+E)^{-1} - A^{-1}\|}{\|A^{-1}\|} \leq
\|A^{-1}\|\|E\|+O(\|E\|^2).
\]
Thus, a condition measure for this problem may be taken as
$\|A^{-1}\|$. Associated with this is the classical Eckart-Young
theorem found in \cite{Eckart}, relating the above condition measure
to the distance to ill-posedness.

\begin{thm}[Eckart-Young]\label{EY} For any non-singular matrix, $A$,
\[
\min_{E}\{ \|E\| : A + E \mbox{ is singular}\} ~= ~\frac{1}{\|A^{-1}\|}.
\]
\end{thm}

We are typically concerned with relative condition numbers as
introduced by Demmel in \cite{Demmel2}. For example, with respect to
the problem of matrix inversion, the relative condition number is
$k(A) := \|A\| \|A^{-1}\|$, the commonly used condition measure.

Condition numbers are also important from an algorithmic
perspective. In the above example of matrix inversion, for instance,
the sensitivity of a problem under perturbations could come into
prominence regarding errors in either the initial problem data or
accumulated computational error due to rounding. Hence, it seems
natural that condition numbers should affect algorithm speed. For
example, in the context of linear programming, Renegar defined a
condition measure based on the distance to ill-posedness in
\cite{Renegar1}--in a similar sense as the Eckart-Young result--and
showed its effect on the convergence rate of interior point methods
in \cite{Renegar2}.

For another example, consider the problem of finding a solution to
the system $Ax = b$ where $A$ is a positive-definite matrix. It was
shown in \cite{Akaike} that the steepest descent method is linearly
convergent with rate $(\frac{k(A)-1}{k(A)+1})^2$ and that this bound
is asymptotically tight for almost all choices of initial iterates.
Similarly, it is well known (see \cite{Golub}) that the conjugate
gradient method applied to the same problem is also linearly
convergent with rate $\frac{\sqrt{k(A)} - 1}{\sqrt{k(A)}+1}$.

From a computational perspective, a related and important area of
study is that of error bounds. Given a subset of a Euclidean space,
an error bound is an inequality that bounds the distance from a test
vector to the specified subset in terms of some residual function
that is typically easy to compute. In that sense, an error bound can
be used both as part of a stopping rule during implementation of an
algorithm as well as an aide in proving algorithmic convergence. A
comprehensive survey of error bounds for a variety of problems
arising in optimization can be found in \cite{Pang}.

With regards to the problem of solving a nonsingular linear system
$Ax = b$, one connection between condition measures and error bounds
is immediate. Let $x^*$ be a solution to the system and $x$ be any
other vector. Then
\[
\|x-x^*\| = \|A^{-1}A(x-x^*)\| = \|A^{-1}(Ax-b)\|\leq
\|A^{-1}\|\|Ax-b\|,
\]
so the distance to the solution set is bounded by a constant
multiple of the residual vector, $\|Ax-b\|$, and this constant is
the same one that appears in the context of conditioning and
distance to infeasibility. As we discuss later, this result is not confined to systems of linear equations.

As a result, error bounds and the related condition numbers often
make a prominent appearance in convergence proofs for a variety of
algorithms. In this paper, motivated by a recent randomized iterated projection scheme for systems of linear equations due to Strohmer and Vershynin in
\cite{Strohmer},  we revisit some classical algorithms and
show that, with an appropriate randomization scheme, we can
demonstrate convergence rates directly in terms of these natural
condition measures. The rest of the paper is organized as follows. In the
remainder of this section, we define some notation used
throughout the rest of this paper. In Section \ref{coordinate}, we consider the
problem of solving a linear system $Ax = b$ and show that a
randomized coordinate descent scheme, implemented according to a
specific probability distribution, is linearly convergent with a
rate expressible in terms of traditional conditioning measures. In
Section \ref{alternating}, we build upon the work of Strohmer and Vershynin in
\cite{Strohmer} by considering randomized iterated projection
algorithms for linear inequality systems. In particular, we show how
randomization can provide convergence rates in terms of the
traditional Hoffman error bound in \cite{Hoffman} as well as in
terms of Renegar's distance to infeasibility from \cite{Renegar3}.
In Section \ref{metric}, we consider randomized iterated projection
algorithms for general convex sets and, under appropriate metric
regularity assumptions, obtain {\em local} convergence rates in terms of the
modulus of regularity.

The classical, deterministic versions of the simple algorithms we
consider here have been widely studied, in part due to the extreme simplicity of each iteration:  their linear convergence is well-known.  However, as remarked for linear systems of equations in \cite{Strohmer}, randomized versions are interesting for several
reasons.  The randomized iterated projection method for linear equations from which this work originated may have some practical promise, even compared with conjugate gradients, for example \cite{Strohmer}.  Our emphasis here, however, is theoretical:  randomization here provides a framework for simplifying the
analysis of algorithms, allowing easy bounds on the rates of
linear convergence in terms of natural linear-algebraic condition measures, such as
relative condition numbers, Hoffman constants, and the modulus of
metric regularity.

\section{Notation}
On the Euclidean space $\Rn$, we denote the Euclidean norm by $\|\cdot\|$.  Let $e_i$ denote the column vector with a 1 in the $i^{th}$ position and zeros elsewhere.

We consider $m$-by-$n$ real matrices $A$.  We denote the set of rows of $A$ by $\{a_1^T,\ldots,a_m^T\}$
and the set of columns is denoted $\{A_1,\ldots,A_n\}$.
The {\em spectral norm} of $A$ is the quantity
$\|A\|_2 := \max_{\|x\|=1}\|Ax\|$ and the {\em Frobenius norm} is
$\|A\|_F := \sum_{i,j} a_{ij}^2$.  These norms satisfy the following inequality \cite{Horn}:
\bmye \label{upper}
\|A\|_F \le \sqrt{n} \|A\|_2.
\emye

For an arbitrary matrix, $A$, let $\|A^{-1}\|_2$ be the smallest
constant $M$ such that $\|Ax\|_2 \geq \frac{1}{M}\|x\|_2$ for all
vectors $x$.  In the case $m \ge n$, if $A$ has singular values
$\sigma_1 \ge \sigma_2 \ge \cdots \ge \sigma_n$, then
$M$ can also be expressed as the reciprocal
of the minimum singular value $\sigma_n$, and, if $A$ is invertible, this
quantity equals the spectral norm of $A^{-1}$.

The {\em relative condition number} of $A$ is the quantity
$k(A) := \|A\|_2\|A^{-1}\|_2$; related to this is the {\em scaled condition number}
introduced by Demmel in \cite{Demmel}, defined by
$\kappa(A) := \|A\|_F\|A^{-1}\|_2$. From this, it is easy to verify (using the singular value decomposition, for example) the following relationship between condition numbers:
\[
1 \leq \frac{\kappa(A)}{\sqrt{n}} \leq k(A).
\]

Now suppose the matrix $A$ is $n$-by-$n$ symmetric and positive
definite. The {\em energy norm} (or $A$-norm), denoted
$\|\cdot\|_A$, is defined by $\|x\|_A := \sqrt{x^TAx}$.  The
inequality
\bmye \label{energy}
\|x\|_A^2 \le \|A^{-1}\|_2 \cdot \|Ax\|^2 ~~\mbox{for all}~ x \in \Rn
\emye
is useful later. Further, if $A$ is simply positive semi-definite,
we can generalize Inequality \ref{energy}:
\bmye \label{energy2}
x^TAx \leq \frac{1}{\underline{\lambda}(A)} \|Ax\|^2
\emye
where $\underline{\lambda}(A)$ is the smallest non-zero eigenvalue
of $A$.  We denote the trace of $A$ by $\tra A$:  it satisfies the inequality
\bmye \label{lower} \|A\|_F \ge \frac{\tra A}{\sqrt{n}}. \emye

Given a nonempty closed convex set $S$, let $P_S(x)$ be the projection of
$x$ onto $S$:  that is, $P_S(x)$ is the vector $y$ that is the optimal
solution to $\min_{z\in S} \|x-z\|_2$. Additionally, define the
distance from $x$ to a set $S$ by
\[
d(x,S) = \min_{z\in S} \|x-z\|_2 = \|x-P_S(x)\|.
\]
The following useful inequality is standard:
\bmye \label{pythagoras}
\|y - x\|^2 - \|P_S(y) - x\|^2 \ge  \|y - P_S(y)\|^2 ~~ \mbox{for all}~ x \in S,~ y \in \Rn.
\emye

\section{Randomized Coordinate Descent} \label{coordinate}
Let $A$ be an $n$-by-$n$ positive-definite matrix.  We consider a linear system of the form $Ax=b$, with solution $x^* = A^{-1}b$.  We consider the
equivalent problem of minimizing the strictly convex quadratic
function
\[
f(x) = \frac{1}{2}x^TAx - b^Tx,
\]
and note the standard relationship
\bmye \label{standard}
f(x) - f(x^*) =  \frac{1}{2} \|x-x^*\|_A^2.
\emye

Suppose our current
iterate is $x$ and we obtain a new iterate $x_+$ by performing an
exact line search in the nonzero direction $d$:  that is, $x_+$ is the solution
to $\min_{t\in\R} f(x+td)$. This gives us
\[
x_+ = x + \frac{(b-Ax)^Td}{d^TAd}d
\]
and
\bmye \label{FUpd1}
f(x_+) - f(x^*) = \frac{1}{2}\|x_+-x^*\|^2_A =
\frac{1}{2}\|x-x^*\|^2_A - \frac{((Ax-b)^Td)^2}{2d^TAd}.
\emye
One natural
choice of a set of easily-computable search directions is to choose
$d$ from the set of coordinate directions, $\{e_1,\ldots,e_n\}$.
Note that, when using search direction $e_i$, we can compute the new
point
\[
x_+ = x + \frac{b_i-a_i^Tx}{a_{ii}}e_i
\]
using only $2n + 2$ arithmetic operations. If the search direction
is chosen at each iteration by successively cycling through the set
of coordinate directions, then the algorithm is known to be linearly
convergent but with a rate not easily expressible in terms of
typical matrix quantities (see \cite{Golub} or \cite{Quarteroni}) .
However, by choosing a coordinate direction as a search direction
randomly according to an appropriate probability distribution, we
can obtain a convergence rate in terms of the relative condition
number. This is expressed in the following result.

\begin{alg}\label{CDPD}
Consider an $n$-by-$n$ positive semidefinite system $Ax=b$ and let $x_0 \in \Rn$ be an arbitrary starting point. For $j=0,1,\ldots,$ compute
\[
x_{j+1} = x_j + \frac{b_i - a_i^Tx_j}{a_{ii}}e_i
\]
where, at each iteration $j$, the index $i$ is chosen independently at random from the set $\{1,\ldots,n\}$, with distribution
\[
P\{ i = k \} ~=~ \frac{a_{kk}}{\tra A}.
\]
\end{alg}

Notice in the algorithm that the matrix $A$ may be singular, but that nonetheless $a_{ii} > 0$ almost surely.  If $A$ is merely positive semidefinite, solutions of the system $Ax=b$ coincide with minimizers of the function $f$, and consistency of the system is equivalent to $f$ being bounded below.   We now have the following result.

\begin{thm}\label{LC1}
Consider a consistent positive-semidefinite system $Ax=b$, and define the corresponding objective and error by
\begin{eqnarray*}
f(x) & = & \frac{1}{2}x^TAx - b^Tx \\
\delta(x) & = & f(x) - \min f.
\end{eqnarray*}
Then Algorithm \ref{CDPD} is linearly convergent in expectation:  indeed,
for each iteration $j=0,1,2,\ldots$,
\[
E[\delta(x_{j+1}) \:|\: x_j] ~\leq~
\Big(1-\frac{\underline{\lambda}(A)}{\tra A} \Big) \delta(x_j).
\]
In particular, if $A$ is positive-definite and $x^* = A^{-1}b$, we have the equivalent property
\[
E[\|x_{j+1}-x^*\|^2_{A} \:|\: x_j] ~\leq~
\Big( 1-\frac{1}{\|A^{-1}\|_2 \tra A} \Big) \|x_j-x^*\|^2_{A}.
\]
Hence the expected reduction in the squared error
$\|x_j - x^*\|_A^2$ is at least a factor
\[
1-\frac{1}{\sqrt{n}\kappa(A)} ~\le~ 1-\frac{1}{n k(A)}
\]
at each iteration.
\end{thm}

\pf
Note that if coordinate direction
$e_i$ is chosen during iteration $j$, then Equation \ref{FUpd1}
shows
\[
f(x_{j+1}) = f(x_j) - \frac{(b_i - a_i^Tx_j)^2}{2a_{ii}}.
\]
Hence, it follows that
\[
E[f(x_{j+1}) ~|~ x_j] = f(x_j) - \sum_{i=1}^n
\frac{a_{ii}}{\tra(A)}\frac{(b_i - a_i^Tx_j)^2}{2a_{ii}} = f(x_j)
- \frac{1}{2\tra A}\|Ax_j-b\|^2.
\]
Using Inequality \ref{energy2} and Equation \ref{standard}, we
easily verify
\[
\frac{1}{2}\|Ax_j-b\|^2 \ge \underline{\lambda}(A) \delta(x_j),
\]
and the first result follows.  Applying Equation
\ref{standard} provides the second result. The final result comes
from applying Inequalities \ref{upper} and \ref{lower}.
\finpf

The simple idea behind the proof of Theorem \ref{LC1} is the main engine driving the
remaining results in this paper. Fundamentally, the idea is to
choose a probability distribution so that the expected distance to
the solution from the new iterate is the distance to the solution
from the old iterate minus some multiple of a residual. Then, using
some type of error bound to bound the distance to a solution in
terms of the residual, we obtain expected linear convergence of the
algorithm.

Now let us consider the more general problem of finding a solution
to a linear system $Ax = b$ where $A$ is an $m \times n$.  More generally, since the system might be inconsistent, we seek a ``least squares solution'' by minimizing the function
$\|Ax-b\|^2$.  The minimizers are exactly the solutions of the positive-semidefinite system $A^TAx = A^Tb$, to which we could easily apply the previous algorithm; however, as usual, we wish to avoid computing the new matrix $A^TA$ explicitly. Instead, we can proceed as follows.

\begin{alg}\label{CDGen}
Consider a linear system $Ax = b$ for a nonzero $m$-by-$n$ matrix $A$. Let $x_0 \in \Rn$ be an arbitrary initial point and let $r_0 = b - Ax_0$ be the initial residual. For each $j=0,1,\ldots,$ compute
\begin{eqnarray*}
\alpha_j & = & \frac{A_i^Tr_j}{\|A_i\|^2} \\
x_{j+1} & = & x_j + \alpha_j e_i \\
r_{j+1} & = & r_j - \alpha_j A_i,
\end{eqnarray*}
where, at each iteration $j$, the index $i$ is chosen independently at random from the set $\{1,\ldots,n\}$, with distribution
\[
P\{ i = k \} ~=~ \frac{\|A_k\|^2}{\|A\|_F^2} ~~~ (k = 1,2,\ldots,n).
\]
\end{alg}

\noindent
(In the formula for $\alpha_j$, notice by assumption that $A_i \ne 0$ almost surely.)

Note that the step size at each iteration can be obtained by
directly minimizing the residual in the respective coordinate
direction. However, the algorithm can also be viewed as the
application of the algorithm for positive definite systems on the
system of normal equations, $A^TAx = A^Tb$, without actually having
to compute the matrix $A^TA$. Given the motivation of directly
minimizing the residual, we would expect that Algorithm
\ref{CDGen} would converge to a least squares solution, even in the case where the underlying system is inconsistent. The next result shows that this
is, in fact, the case.

\begin{thm}\label{CDMN2}
Consider any linear system $Ax=b$, where the matrix $A$ is nonzero.  Define the least-squares residual and the error by
\begin{eqnarray*}
f(x) & = & \frac{1}{2} \|Ax-b\|^2 \\
\delta(x) & = & f(x) - \min f.
\end{eqnarray*}
Then Algorithm \ref{CDGen}~is linearly convergent in
expectation to a least squares solution for the system: for each
iteration $j=0,1,2\ldots$,
\[
E[\delta(x_{j+1}) \:|\: x_j] ~\leq~
\Big(1-\frac{\underline{\lambda}(A^TA)}{\|A\|^2_F}\Big) \delta(x_j).
\]
In particular, if $A$ has full column rank, we have the equivalent property
\[
E[\|x_{j+1}-\hat{x}\|^2_{A^TA} \:|\: x_j] ~\leq~
\Big(1-\frac{1}{\kappa(A)^2}\Big)\|x_j-\hat{x}\|^2_{A^TA}
\]
where $\hat{x} = (A^TA)^{-1}A^Tb$ is the unique least-squares
solution.
\end{thm}

\pf It is easy to verify, by induction on $j$, that the iterates
$x_j$ are exactly the same as the iterates generated by Algorithm
\ref{CDPD}, when applied to the positive semi-definite system $A^TAx
= A^Tb$, and furthermore that the residuals satisfy $r_j = b-Ax_j$
for all $j=0,1,2,\ldots$.  Hence, the results follow directly by
Theorem \ref{LC1}.
\finpf

By the coordinate descent nature of this algorithm, once we have
computed the initial residual $r_0$ and column norms
$\{\|A_i\|^2\}_{i=1}^n$, we can perform each iteration in $O(n)$
time, just as in the positive-definite case. Specifically, this new
iteration takes $4n+1$ arithmetic operations, compared with $2n+2$
for the positive-definite case.

For a computational example, we apply Algorithm \ref{CDGen} to
random $500\times n$ matrices where each element of $A$ and $b$ is
an independent Gaussian random variable and we let $n$ take values
50, 100, 150 and 200.

\begin{center}
\includegraphics[width=15.5cm]{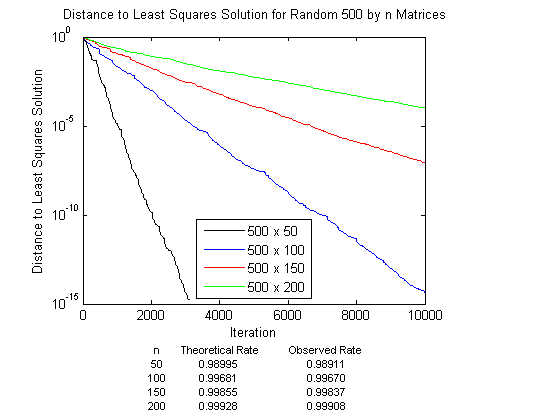}
\end{center}

\noindent Note that in the above examples, the theoretical bound
provided by Theorem \ref{CDMN2} predicts the actual behavior of the
algorithm reasonably well.

\section{Randomized Iterated Projections} \label{alternating}
Iterated projection algorithms share some important characteristics
with coordinate descent algorithms. Both are well studied and much
convergence theory exists; a comprehensive overview on iterated
projections can be found in \cite{Deutsch}. However, even for linear
systems of equations, standard developments do not provide bounds on
convergence rates in terms of usual condition numbers. By contrast,
in the recent paper \cite{Strohmer}, Strohmer and Vershynin obtained
such bounds via the following randomized iterated projection
algorithm, which also provided the motivation for our work in the
previous section.

\begin{alg} \label{SV1}
Consider a linear system $Ax = b$ for a nonzero $m$-by-$n$ matrix $A$. Let $x_0 \in \Rn$ be an arbitrary initial point. For each $j=0,1,\ldots,$ compute
\[
x_{j+1} = x_j - \frac{a_i^Tx_j - b_i}{\|a_i\|^2}a_i
\]
where, at each iteration $j$, the index $i$ is chosen independently at random from the set $\{1,\ldots,m\}$, with distribution
\[
P\{ i = k \} ~=~ \frac{\|a_k\|^2}{\|A\|_F^2} ~~~ (k = 1,2,\ldots,m).
\]
\end{alg}

Notice that the new iterate $x_{j+1}$ is simply the orthogonal projection of the old iterate $x_j$ onto the hyperplane $\{ x : a_i^T x = b_i \}$.  At first sight, the choice of probability distribution may seem curious, since we could rescale the equations arbitrarily without having any impact on the projection operations.  However, following \cite{Strohmer}, we emphasize that the aim is to understand linear convergence rates in terms of {\em linear-algebraic} condition measures associated with the original system, rather than in terms of {\em geometric} notions associated with the hyperplanes.  This randomized algorithm has the following behavior.

\begin{thm}[Strohmer-Vershynin, \cite{Strohmer}] \label{AP}
Given any matrix $A$ with full column
rank, suppose the linear system $Ax=b$ has solution $x^*$.  Then
Algorithm \ref{SV1}~converges linearly in expectation:  for each iteration $j=0,1,2,\ldots$,
\[
E[\|x_{j+1}-x^*\|_2^2 ~|~x_j] ~\leq~
\Big(1-\frac{1}{\kappa(A)^2}\Big)\|x_j-x^*\|^2_2
\]
\end{thm}

We seek a way of generalizing the above algorithm and
convergence result to more general systems of linear inequalities, of the form
\bmye \label{system}
\left\{
\begin{array}{rcll}
a_i^T x & \le & b_i & (i \in I_{\le}) \\
a_i^T x &  =  & b_i & (i \in I_{=}),
\end{array}
\right.
\emye
where the disjoint index sets $I_{\le}$ and $I_=$ partition the set $\{1,2,\ldots,m\}$.  To do so, staying with the techniques of the previous section, we need a corresponding error bound for a system of linear inequalities. First, given a vector $x\in\Rn$, define the vector
$x^+$ by $(x^+)_i = \max\{x_i,0\}$. Then a starting point for this
subject is a result by Hoffman in \cite{Hoffman}.

\begin{thm}[Hoffman]\label{HoffC}
For any right-hand side vector $b \in \R^m$, let $S_b$ be the set of feasible solutions of the linear system (\ref{system}).  Then there exists a constant $L$, independent of $b$, with the following property:
\bmye \label{hoffman}
x \in \Rn ~\mbox{and}~ S_b \ne \emptyset ~~\Rightarrow~~ d(x,S_b) \le L \|e(Ax-b)\|,
\emye
where the function $e \colon \R^m \to \R^m$ is defined by
\[
e(y)_i =
\left\{
\begin{array}{cl}
y_i^+ & (i \in I_{\le}) \\
y_i   & (i \in I_{=}),
\end{array}
\right.
\]
\end{thm}

In the above result, each component of the vector $e(Ax-b)$ indicates the error in the corresponding inequality or equation.  In particular $e(Ax-b)=0$ if and only if $x \in S_b$.  Thus Hoffman's result provides a linear bound for the distance from a trial point $x$ to the feasible region in terms of the size of the ``a posteriori error'' associated with $x$.

We call the minimum constant $L$ such that property (\ref{hoffman}) holds the {\em Hoffman constant} for the system (\ref{system}).  Several authors give geometric or algebraic meaning to this constant, or exact expressions for it, including \cite{Guler},
\cite{Ng}, \cite{Li}; for a more thorough treatment of the subject,
see \cite{Pang}.  In the case of linear equations (that is,
$I_{\le} = \emptyset$), an easy calculation using the singular value decomposition shows that the Hoffman constant is just the reciprocal of the smallest nonzero singular value of the matrix $A$, and hence equals $\|A^{-1}\|_2$ when $A$ has full column rank.

For the problem of finding a solution to a system of linear
inequalities, we consider a randomized algorithm generalizing Algorithm \ref{SV1}.

\begin{alg} \label{SV2}
Consider the system of inequalities (\ref{system}).  Let $x_0$ be an arbitrary initial point.
For each $j=0,1,\ldots,$ compute
\begin{eqnarray*}
\beta_j & = &
\left\{
\begin{array}{cl}
(a_i^Tx_j - b_i)^+ & (i \in I_{\le}) \\
a_i^Tx_j - b_i & (i \in I_=)
\end{array}
\right.
\\
x_{j+1} & = & x_j - \frac{\beta_j}{\|a_i\|^2}a_i
\end{eqnarray*}
where, at each iteration $j$, the index $i$ is chosen independently at random from the set $\{1,\ldots,m\}$, with distribution
\[
P\{ i = k \} ~=~ \frac{\|a_k\|^2}{\|A\|_F^2} ~~~ (k = 1,2,\ldots,m).
\]
\end{alg}

\noindent
In the above algorithm, notice $\beta_j = e(Ax_j-b)_i$.  We can now generalize Theorem \ref{AP} as follows.

\begin{thm}\label{APHoff}
Suppose the system (\ref{system}) has nonempty feasible region $S$.  Then Algorithm \ref{SV2} converges linearly in
expectation: for each iteration $j=0,1,2,\ldots$,
\[
E[d(x_{j+1},S)^2 ~|~ x_j] ~\leq~ \Big(1 - \frac{1}{L^2\|A\|^2_F}
\Big) d(x_j,S)^2.
\]
where $L$ is the Hoffman constant.
\end{thm}

\pf Note that if the index $i$ is chosen during iteration $j$, then it
follows that
\begin{eqnarray*}
\lefteqn{
\|x_{j+1}-P_S(x_{j+1})\|^2_2  ~ \leq ~ \|x_{j+1}-P_S(x_j)\|^2_2
}  \\
& & = ~
\Big\| x_j - \frac{e(Ax_j-b)_i}{\|a_i\|^2}a_i -P_S(x_j) \Big\|^2_2 \\
& & = ~
\| x_j - P_S(x_j) \|^2_2 + \frac{e(Ax_j-b)_i^2}{\|a_i\|^2}
- 2 \frac{e(Ax_j-b)_i}{\|a_i\|^2}a_i^T(x_j - P_S(x_j)).
\end{eqnarray*}
Note $P_S(x_j)\in S$.  Hence if $i \in I_{\le}$, then
$a_i^T P_S(x_j)\leq b_i$, and $e(Ax_j-b)_i \ge 0$, so
\[
e(Ax_j-b)_i a_i^T(x_j - P_S(x_j)) \ge e(Ax_j-b)_i (a_i^T x_j - b_i) = e(Ax_j-b)_i^2.
\]
On the other hand, if $i \in I_=$, then
$a_i^T P_S(x_j) = b_i$, so
\[
e(Ax_j-b)_i a_i^T(x_j - P_S(x_j)) = e(Ax_j-b)_i (a_i^T x_j - b_i) = e(Ax_j-b)_i^2.
\]
Putting these two cases together with the previous inequality shows
\[
d(x_{j+1},S)^2 ~\le~ d(x_j,S)^2 - \frac{e(Ax_j-b)_i^2}{\|a_i\|^2}.
\]
Taking the expectation with respect to the specified probability
distribution, it follows that
\[
E[d(x_{j+1},S)^2 ~|~ x_j] \le d(x_j,S)^2 -
\frac{\|e(Ax_j-b)\|^2}{\|A\|_F^2}
\]
and the result now follows by the Hoffman bound. \finpf

\noindent
Since Hoffman's bound is not independent of the scaling of the
matrix $A$, it is not surprising that a normalizing constant like $\|A\|^2_F$ term appears in the result.

For a computational example, we consider linear inequality systems
$Ax \leq b$ where the elements of $A$ are independent standard
Gaussian random variables and $b$ is chosen so that the resulting
system has a non-empty interior. We consider matrices $A$ which are
$500\times n$, letting $n$ take values 50, 100, 150 and 200. We then
apply Algorithm \ref{SV2} to these problems and observe the
following computational results.

\begin{center}
\includegraphics[width=13.5cm]{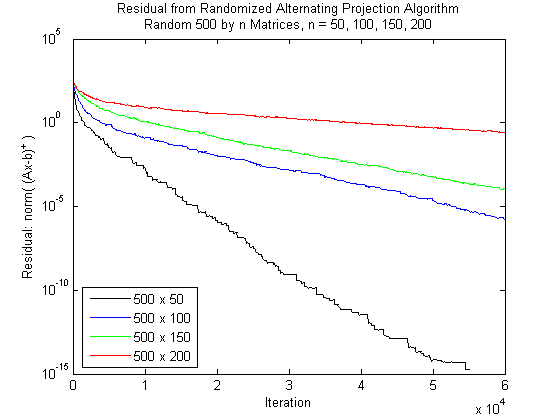}
\end{center}

Another natural conditioning measure for linear inequality systems
is the distance to infeasibility, defined by Renegar in
\cite{Renegar3}, and shown in \cite{Renegar2} to govern the
convergence rate of interior point methods for linear programming.
It is interesting, therefore, from a theoretical perspective, to
obtain a linear convergence rate for iterated projection
algorithms in terms of this condition measure as well. For simplicity, we concentrate on the inequality case, $Ax \le b$.  To begin, let
us recall the following results.

The {\em distance to infeasibility} \cite{Renegar3} for the system
$Ax \le b$ is the number
\[
\mu ~=~ \inf \Big\{ \max \{ \| \Delta A \| , \| \Delta b \| \} : (A
+ \Delta A) x \le b + \Delta b ~ \mbox{is infeasible} \Big\}.
\]

\begin{thm}[Renegar, \cite{Renegar3}, Thm 1.1] \label{renegar}
Suppose $\mu > 0$.  Then there exists a point $\hat x \in S$
satisfying $\|\hat x\| \le \|b\|/\mu$.  Furthermore, any point $x
\in \R^n$ satisfies the inequality
\[
d(x,S) ~\le~ \frac{\max \{ 1 , \|x\| \}}{\mu} \|(Ax-b)^+\|.
\]
\end{thm}

Using this, we can bound the linear convergence rate for the
Algorithm \ref{SV2} in terms of the distance to infeasibility, as follows.
Notice first that $\|x_j - \hat x\|$ is nonincreasing in $j$, by Inequality \ref{pythagoras}.  Suppose we start Algorithm \ref{SV2} at the initial point $x_0 = 0$. Applying Theorem \ref{renegar}, we see that for all $j=1,2,\ldots,$
\[
\|x_j\| \le \|\hat x\| + \|x_j  - \hat x\| \le
\|\hat x\| + \|x_0 - \hat x\| \le \frac{2\|b\|}{\mu},
\]
so
\[
d(x_j , S) \le \max \Big\{ \frac{1}{\mu} , \frac{2\|b\|}{\mu^2}
\Big\} \|(Ax_j-b)^+\|.
\]
Using this inequality in place of Hoffman's bound in the proof of
Theorem \ref{APHoff} gives
\[
E[d(x_{j+1},S)^2 ~|~ x_j] \leq \left[1 -
\frac{1}{\|A\|_F^2(\max\{\frac{1}{\mu},
\frac{2\|b\|}{\mu^2}\})^2}\right]d(x_j,S)^2.
\]
Although this bound may not be the best possible (and, in fact, it
may not be as good as the bound provided in Theorem \ref{APHoff}),
this result simply emphasizes a relationship between algorithm speed
and conditioning measures that appears naturally in other contexts.
In the next section, we proceed with these ideas in a more general
framework.

\section{Metric Regularity and Local Convergence} \label{metric}
The previous section concerned global rates of linear convergence.  If instead we are interested in {\em local} rates, we can re-examine a generalization of our problem
through an alternative perspective of set-valued mappings. Consider
a set-valued mapping $\Phi: \Rn\tto \Rm$ and the problem of solving
the associated constraint system of the form $b\in \Phi(x)$ for the
unknown vector $x$. For example, finding a feasible solution to
$Ax\leq b$ is equivalent to finding an $x$ such that
\bmye  \label{inclusion}
b \in Ax + \Rm_+.
\emye

Related to this is the idea of \textit{metric regularity} of
set-valued mappings. We say the set-valued mapping $\Phi$ is
metrically regular at $\bar x$ for $\bar b\in \Phi(\bar x)$ if there
exists $\gamma > 0$ such that \bmye \label{MetReg} d(x,
\Phi^{-1}(b)) \leq \gamma d(b, \Phi(x)) \mbox{ for all } (x,b)
\mbox{ near } (\bar x, \bar b), \emye where $\Phi^{-1}(b) = \{x:
b\in \Phi(x)\}$. Further, the \textit{modulus of regularity} is the
infimum of all constants $\gamma$ such that Equation \ref{MetReg}
holds.  Metric regularity is strongly connected with a variety of ideas from
variational analysis: a good background reference is \cite{Roc98}.

Metric regularity generalizes the error bounds discussed in previous
sections at the expense of only guaranteeing a bound in local terms.
For example, if $\Phi$ is a single-valued linear map, then the
modulus of regularity (at any $\bar{x}$ for any $\bar{b}$)
corresponds to the typical conditioning measure $\|\Phi^{-1}\|$
(with $\|\Phi^{-1}\| = \infty$ implying the map is not metrically
regular) and if $\Phi$ is a smooth single-valued mapping, then the
modulus of regularity is the reciprocal of the minimum singular
value of the Jacobian, $\nabla \Phi(x)$. From an alternative
perspective, metric regularity provides a framework for generalizing
the Eckart-Young result on the distance to ill-posedness of linear
mappings cited in Theorem \ref{EY}. Specifically, if we
define the \textit{radius of metric regularity} at $\bar{x}$ for
$\bar{b}$ for a set-valued mapping $\Phi$ between finite dimensional
spaces by
\[
\mbox{rad}\Phi(\bar{x}|\bar{b}) ~=~
\inf \{\|E\| : \Phi + E \mbox{
not metrically regular at } \bar{x} \mbox{ for } \bar{b} +
E(\bar{x}) \},
\]
where the infimum is over all linear functions $E$,
then one obtains the strikingly simple relationship (see \cite{Dontchev})
\[
\mbox{modulus of regularity of $\Phi$  at $\bar{x}$ for $\bar{b}$}
~=~ \frac{1}{\mbox{rad}\Phi(\bar{x}|\bar{b})}.
\]
We will not be directly using the above result. Here, we simply use the fundamental idea of metric regularity which says that the distance from a point to the solution set,
$d(x,\Phi^{-1}(b))$, is locally bounded by some constant times a
''residual''. For example, in the case where $\Phi$ corresponds to
the linear inequality system (\ref{inclusion}), we have that
$d(b,\Phi(x)) = \|(Ax-b)^+\|$ implies that the modulus of regularity
is in fact a global bound and equals the Hoffman bound. More
generally, we wish to emphasize that metric regularity ties together
several of the ideas from previous sections at the expense of those
results now only holding locally instead of globally.

In what follows, assume all distances are Euclidean distances. We
wish to consider how the modulus of regularity of $\Phi$ affects the
convergence rate of iterated projection algorithms.  We remark that
linear convergence for iterated projection methods on convex sets
has been very widely studied:  see \cite{Deutsch}, for example. Our
aim here is to observe, by analogy with previous sections, how
randomization makes the linear convergence rate easy to interpret in
terms of metric regularity.

Let
$S_1,S_2,\ldots,S_m$ be closed convex sets in a Euclidean space $\E$
such that $\cap_i S_i\neq \emptyset$. Then, in a manner similar to
\cite{Lewis}, we can endow the product space
$\E^m$ with the inner product
\[
\langle (u_1,u_2,\ldots,u_m), (v_1,v_2,\ldots, v_m)\rangle  =
\sum_{i=1}^m \langle u_i,v_i\rangle
\]
and consider the set-valued mapping $\Phi: \E\rightarrow\E^m$ given
by
\bmye \label{PhiMap}
\Phi(x) = (S_1 - x, S_2 - x, \ldots, S_m -x).
\emye
Then it clearly follows that $\bar{x}\in\cap_i S_i \Leftrightarrow
0\in\Phi(\bar{x})$. Under appropriate regularity assumptions, we
obtain the following local convergence result.

\begin{thm}\label{MRAP}
Suppose the set-valued mapping $\Phi$ given by Equation \ref{PhiMap}
is metrically regular at $\bar{x}$ for 0 with regularity modulus
$\gamma$. Let $\bar{\gamma}$ be any constant strictly larger than
$\gamma$ and let $x_0$ be any initial point sufficiently close to
$\bar{x}$. Further, suppose that $x_{j+1} = P_{S_i}(x_j)$ with
probability $\frac{1}{m}$ for $i=1,\ldots,m$. Then
\[
E[d(x_{j+1}, S)^2 ~|~ x_j] ~\leq~
\Big(1-\frac{1}{m\bar{\gamma}^2}\Big)d(x_j,S)^2.
\]
\end{thm}

\pf
First, note that by Inequality \ref{pythagoras}, the distance $\|x_j - \bar x\|$ is nonincreasing in $j$.  Hence if $x_0$ is
sufficiently close to $\bar{x}$, then $x_j$ is as well for all
$j\geq 0$. Then, again using Inequality \ref{pythagoras} (applied to the set $S_i$), we have, for all points $x \in S \subset S_i$,
\[
\|x_j - x\|^2 - \|x_j - P_{S_i}(x_j)\|^2 \ge \|P_{S_i}(x_j) - x\|^2.
\]
Taking the minimum over $x \in S$, we deduce
\[
d(x_j,S)^2 - \|x_j - P_{S_i}(x_j)\|^2 \ge d(P_{S_i}(x_j),S)^2.
\]
Hence
\begin{eqnarray*}
E[d(x_{j+1}, S)^2 ~|~ x_j] & = &
\frac{1}{m}\sum_{i=1}^m d(P_{S_i}(x_j),S)^2 \\
& \leq &
\frac{1}{m}\sum_{i=1}^m [d(x_j,S)^2 - d(x_j, S_i)^2] \\
& = &
d(x_j,S)^2 - \frac{1}{m}\sum_{i=1}^m d(x_j,S_i)^2 \\
& = &
d(x_j,S)^2 - \frac{1}{m}d(0, \Phi(x_j))^2 \\
& \leq &
\Big(1-\frac{1}{m\bar{\gamma}^2}\Big)d(x_j,S)^2,
\end{eqnarray*}
using the definition of metric regularity.
\finpf

Note that metric regularity at $\bar{x}$ for 0 is a slightly
stronger assumption than actually necessary for this result.
Specifically, the above result holds as long as Equation
\ref{MetReg} holds for all $x$ near $\bar{x}$ with $\bar{b}=0$
fixed, as opposed to the above definition requiring it to hold for
all $b$ near $\bar{b}$ as well.

For a moment, let $m=2$ and consider the sequence of iterates
$\{x_j\}_{j\geq 0}$ generated by the randomized iterated
projection algorithm. By idempotency of the projection operator,
there's no benefit to projecting onto the same set in two
consecutive iterations, so the subsequence consisting of different
iterates corresponds exactly to that of the non-randomized
iterated projection algorithm. In particular, if $x_j\in S_1$,
then
\[
d(P_{S_2}(x_j), S)^2 \leq d(x_j,S)^2 - d(x_j,S_2)^2
= d(x_j,S)^2 - [d(x_j,S_2)^2 + d(x_j,S_1)^2]
\]
since $d(x_j,S_1) = 0$. This gives us the following corollary, which also follows through more standard deterministic arguments.

\begin{cor}\label{MRAPCor}
If $\Phi$ is metrically regular at $\bar{x}$ for 0 with regularity
modulus $\gamma$ and $\bar{\gamma}$ is larger than $\gamma$, then for $x_0$ sufficiently close to $\bar{x}$, the
2-set iterated projection algorithm is linearly convergent and
\[
d(x_{j+1},S)^2 \leq \Big(1-\frac{1}{\bar{\gamma}^2}\Big)d(x_j,S)^2.
\]
\end{cor}

Further, consider the following refined version of the $m$-set randomized
algorithm. Suppose $x_0\in S_1$ and $i_0 = 1$. Then for
$j=1,2,\ldots,$ let $i_j$ be chosen uniformly at random from
$\{1,\ldots,m\}\backslash \{i_{j-1}\}$ and $x_{j+1} =
P_{S_{i_j}}(x_j)$. Then we obtain the following similar result.

\begin{cor}\label{MRAPCor2}
If $\Phi$ is metrically regular at $\bar{x}$ for 0 with regularity modulus
$\gamma$ and $\bar{\gamma}$ is larger than $\gamma$, then for $x_0$ sufficiently close to $\bar{x}$, the refined
$m$-set randomized iterated projection algorithm is linearly
convergent in expectation and
\[
E[d(x_{j+1},S)^2 ~|~ x_j,~ i_{j-1}] \leq
\Big(1-\frac{1}{(m-1)\bar{\gamma}^2}\Big) d(x_j,S)^2.
\]
\end{cor}

A simple but effective product space formulation by Pierra in
\cite{Pierra} has the benefit of reducing the problem of finding a
point in the intersection of finitely many sets to the problem of finding
a point in the intersection of 2 sets. Using the
notation above, we consider the closed set in the
product space given by
\[
T = S_1 \times S_2 \times \ldots \times S_m
\]
and the subspace
\[
L = \{Ax: x\in E\}
\]
where the linear mapping $A:\E\rightarrow\E^m$ is defined by $Ax =
(x,x,\ldots,x)$. Again, notice that $\bar{x}\in \cap_i S_i
\Leftrightarrow (\bar{x},\ldots,\bar{x}) \in T \cap L$. One
interesting aspect of this formulation is that projections in the
product space $\E^m$ relate back to projections in the original
space $\E$ by
\begin{eqnarray*}
(z_1,\ldots,z_m) \in P_T(Ax)
& \Leftrightarrow &
z_i\in P_{S_i}(x) ~~(i=1,2,\ldots,m) \\
(P_L(z_1,\ldots,z_m))_i & = & \frac{1}{m}(z_1+z_2+\ldots+z_m)
~~(i=1,\ldots,m)
\end{eqnarray*}
This formulation provides a nice analytical framework:  we
can use the above equivalence of projections to consider the
\textit{method of averaged projections} directly, defined as
follows.

\begin{alg}\label{AvgP}Let $S_1,\ldots,S_m\subseteq E$ be nonempty closed convex sets. Let $x_0$ be an
initial point. For $j=1,2,\ldots$, let
\[
x_{j+1} = \frac{1}{m}\sum_{i=1}^m P_{S_i}(x_j).
\]
\end{alg}

Simply put, at each iteration, the algorithm projects the current
iterate onto each set individually and takes the average of those
projections as the next iterate. In the product space formulation,
this is equivalent to $x_{j+1} = P_L(P_T(x_j))$. Expanding on the
work of Pierra in \cite{Pierra}, additional convergence theory for
this algorithm has been examined by Bauschke and Borwein in
\cite{Bauschke}. Under appropriate regularity conditions, the
general idea is that convergence of the iterated projection
algorithm for two sets implies convergence of the averaged
projection algorithm for $m$ sets. In a similar sense, we prove the
following result in terms of randomized projections.

\begin{thm}\label{AvgPRandom}
Suppose $S = \cap_{i=1}^m S_i$ is non-empty. If the randomized
projection algorithm of Theorem \ref{MRAP} is linearly convergent in
expectation with rate $\alpha$, then so is Algorithm \ref{AvgP}.
\end{thm}

\pf Let $x_j$ be the current iterate, $x_{j+1}^{AP}$ be the new
iterate in the method of averaged projections and $x_{j+1}^{RP}$ be
the new iterate in the method of uniformly randomized projections.
Then note that:
\[
x_{j+1}^{AP} = \frac{1}{m}\sum_{i=1}^m P_{S_i}(x_j) =
E[x_{j+1}^{RP}].
\]
By convexity of the $S_i$'s, it follows that
\[
d(x_{j+1}^{AP}, S) = d(E[x_{j+1}^{RP}|x_j], S) \leq
E[d(x_{j+1}^{RP}, S) | x_j] \leq \alpha d(x_j, S),
\]
by Jensen's Inequality.
\finpf

Hence, the method of averaged projections converges no more slowly
than the method of uniformly random projections. In particular,
under the assumptions of Theorem \ref{MRAP}, the method of averaged
projections converges with rate no larger than
$1-\frac{1}{m\bar{\gamma}^2}$.

\bibliographystyle{plain}
\bibliography{RandomizedMethods}

\end{document}